\documentclass[12pt]{article}
\usepackage{amsthm,amssymb,amscd,amsmath,latexsym,indentfirst,color,amsfonts}
\usepackage[mathcal]{eucal}
 \usepackage[all]{xy}
\theoremstyle{plain}
\newtheorem{theorem}{Theorem}[section]
\newtheorem{proposition}[theorem]{Proposition}
\newtheorem{coro}[theorem]{Corollary}
\newtheorem{lemma}[theorem]{Lemma}
\newtheorem{ex}[theorem]{Example}

\newtheorem{defi}[theorem]{Definition}

\newcommand{\ble}{\begin {lemma}}
\newcommand{\ele}{\end {lemma}}
\newcommand{\bpr}{\begin {proposition}}
\newcommand{\epr}{\end {proposition}}
\newcommand{\bthm}{\begin {theorem}}
\newcommand{\ethm}{\end {theorem}}
\newcommand{\bco}{\begin {coro}}
\newcommand{\eco}{\end {coro}}
\newcommand{\bex}{\begin {ex}}
\newcommand{\eex}{\end {ex}}

\newcommand{\be}{\begin {equation}}
\newcommand{\ee}{\end {equation}}
\newcommand{\bp}{\begin {proof}}
\newcommand{\ep}{\end {proof}}
\newcommand{\bee}{\begin {equation*}}
\newcommand{\eee}{\end {equation*}}
\newcommand{\rt}{\rightarrow}
\newcommand{\lb}{\label}
\newcommand{\la}{\lambda}

\begin{document}
\title{The Jordan-Chevalley decomposition and Jordan canonical form of a quaternionic linear operator
}
\author{Gang Han\thanks{Corresponding author. The author was supported by
   Zhejiang Province Science Foundation of China, grant No. LY14A010018.}
\\ School of Mathematics \\ Zhejiang University\\mathhgg@zju.edu.cn \\[2mm]
 Jing Yu
\\ School of Mathematics \\ Zhejiang University\\yujing46@126.com\\[2mm] %
Zheyu Sun
\\ School of Mathematics \\ Zhejiang University\\szy4007@outlook.com
}
\date{ }
\maketitle
\begin{abstract}
We introduce some basic notions and results for quaternionic linear operators analogous to those for complex  linear operators. Our main result is to prove the additive and multiplicative Jordan-Chevalley decompositions for quaternionic linear operators, which are related by the exponential map. We also give a new proof of the theorem of Jordan canonical form for quaternionic linear operators, which is intrinsic and takes less computations than  the known proofs for square quaternionic matrices.
\end{abstract}



\section{Introduction}
 \setcounter{equation}{0}\setcounter{theorem}{0}

Quaternions were first described by Irish mathematician William Rowan Hamilton in middle 1900s, and are widely used in both pure and applied mathematics, in particular for calculations involving three-dimensional rotations such as in three-dimensional computer graphics, computer vision, and crystallographic texture analysis. Quaternionic matrices are used in quantum mechanics \cite{ro} and in the treatment of multibody problems \cite{gs}. Most of the algebraic theory on quaternion linear algebra, including quaternionic matrices, can be found in \cite{r}. There is a good survey of some important results on quaternionic matrices in \cite{z1}. The spectral theory for normal quaternionic matrices can be found in \cite{fp}.

The  Jordan canonical form of a linear operator plays a central role in the structure theory of  linear operators on finite-dimensional vector spaces over fields. A further development of Jordan canonical form is the additive Jordan-Chevalley decomposition, which expresses a linear operator as the sum of its commuting semisimple part and nilpotent part. For any linear operator $\phi$ on a finite-dimensional vector space over a perfect field, the Jordan-Chevalley decomposition exists and is unique, i.e. $\phi=\phi_{s}+\phi_n$, where $\phi_{s}$ and $\phi_{n}$ are the respective semisimple part and  nilpotent part, which are in fact expressible as polynomials in $\phi$. A proof of this result in the case the field is algebraically closed can be found in  \cite{h}. There is also a multiplicative Jordan-Chevalley decomposition for  invertible operators on  finite-dimensional vector spaces over  perfect fields.  The decomposition is important in the study of Lie algebras and algebraic groups. We will prove the additive and multiplicative Jordan-Chevalley decompositions for quaternionic linear operators, which are related by the exponential map. This is the main result of the paper.

In \cite{r}, the Jordan canonical form of a square quaternionic matrix is referred to as a key result of quaternion linear algebra, and is proved in Section 5.6. There have been several proofs of the theorem of Jordan canonical form of a square quaternionic matrix \cite{hl}\cite{chw}\cite{jz} \cite{zw}\cite{r}\cite{l}. ( The proof in \cite{zw} is in fact wrong). It took many computations in most of the proofs. The proof in \cite{zw} is relatively simple. We give a new proof of the theorem of Jordan canonical form in terms of quaternionic linear operators, which takes less computations than the known proofs.

In Section 2, we will introduce some basic notions and results of  quaternionic vector spaces and quaternionic linear operators, such as the relations between quaternionic linear operators and complex linear operators, the matrix and the eigenvalues of quaternionic linear operators.
In Section 3, the exponential map for quaternionic linear operators is introduced. The properties of the generalized eigenspaces of the complex linear operator induced from some quaternionic linear operator are studied. We define the characteristic polynomial of a quaternionic linear operator and prove the corresponding Cayley-Hamilton theorem as in \cite{z1}. In Section 4, the additive and multiplicative Jordan-Chevalley decomposition for quaternionic linear operators are proved, and we show that the additive and multiplicative Jordan-Chevalley decomposition of quaternionic linear operators are related by the exponential map. In Section 5, we give a new proof of the theorem of Jordan canonical form for quaternionic linear operators.  As a corollary, the surjectivity of  the exponential map $exp:M_n(\mathbb{H})\rt GL_n(\mathbb{H})$  is proved.

The proofs of some results below for  quaternionic vector spaces and quaternionic linear operators are omitted as it is analogous to those of the corresponding results for complex linear operators.

\section{Some basic notions and results for quaternionic vector spaces and quaternionic linear operators}
 \setcounter{equation}{0}\setcounter{theorem}{0}
Let $\mathbb{R}$ and $\mathbb{C}$ be the field of the real and complex numbers as usual, and $\mathbb{H}$ be the $\mathbb{R}$-algebra of quaternions with basis $1, i, j ,k$.  In the paper, all the vector spaces are finite dimensional, and all the vector spaces over $\mathbb{H}$ or $\mathbb{C}$ will be assumed to be \textbf{right} vector space.

1. First we review some basic results for quaternions. Let \be \eta:\mathbb{C}\rightarrow \mathbb{H}, a+bi\mapsto a+bi, \ee where $a,b\in \mathbb{R}$, be the standard embedding of $\mathbb{C}$ in $\mathbb{H}$, and we will always regard $\mathbb{C}$ as a subalgebra of $\mathbb{H}$.

 Let  \be \varphi:\mathbb{H}\rightarrow \mathbb{C}, q=a+bi+cj+dk\mapsto a+bi， \ee where $a,b,c,d\in \mathbb{R}$. Then $\varphi(q)$ is called the \textit{complex part} of $q$.

For any $q=a+bi+cj+dk\in \mathbb{H}$, $\bar{q}=a-bi-cj-dk$ is the conjugate of $q$; and $|q|=\sqrt{\bar{q}q}= \sqrt{a^2+b^2+c^2+d^2}$ is the norm of $q$.

Two quaternions $p$ and $q$ are \textit{similar} if there exists some nonzero quaternion $u$ such that $u^{-1}pu=q$. For any $q=a_1+b_1i+c_1j+d_1k\in \mathbb{H}$ one knows that the set $[q]$ of quaternions similar to $q$ is \bee [q]=\{a_1+{b_2}i+c_2j+d_2k|b_2^2+c_2^2+d_2^2=b_1^2+c_1^2+d_1^2\}.\eee The complex numbers contained in $[q]$ are $a_1\pm ri,~r=\sqrt{b_1^2+c_1^2+d_1^2}$. See Lemma 2.1 and Theorem 2.2 of \cite{z1}.

The set of quaternions $\mathbb{H}$ can be regarded as a 2-dimensional $\mathbb{C}$-vecotr space, with basis $1,
j$. One has \[\mathbb{H}=1\cdot \mathbb{C}\oplus j\cdot \mathbb{C}.\] It is clear that as a ring, $\mathbb{H}$ is generated by $\mathbb{C}$ and $j$, with
\be\lb{a} zj=j\bar{z},~~ \forall z\in \mathbb{C} ~~ \textrm{and}~~ j^2=-1.\ee

It is clear that the maps $\eta$ and $\varphi$ are $\mathbb{C}$-linear maps. One has $\eta \varphi=I$, $\overline{\eta(z)}=\eta(\bar{z})$ for $z\in \mathbb{C}$, and $\overline{\varphi(q)}={\varphi(\bar{q})}$ for $q\in \mathbb{H}$.\\[3mm]

2. Any $\mathbb{H}$-vector space $V$  can be considered as a $\mathbb{C}$-vector space, denoted by $V_0$, by restricting scalars from $\mathbb{H}$ to $\mathbb{C}$. Any $\mathbb{H}$-linear map $\phi:V\rt W$ between $\mathbb{H}$-vector spaces $V$ and $W$ can be regarded as a $\mathbb{C}$-linear map between $\mathbb{C}$-vector spaces $V_0$ and $W_0$, and denoted by $\phi_0:V_0\rt W_0$.

Let $Vect_\mathbb{H}$ and $Vect_\mathbb{C}$ be the category of right vector spaces over $\mathbb{H}$ and over $\mathbb{C}$ respectively. Then the functor $F$ between $Vect_\mathbb{H}$ and $Vect_\mathbb{C}$ is defined as follows. For any $\mathbb{H}$-vector space $V$ in $Vect_\mathbb{H}$, $F(V)=V_0$. For any $\mathbb{H}$-linear map $\phi:V\rt W$, $F(\phi)=\phi_0:V_0\rt W_0$. It is easy to verify that $F(1_V)=1_{V_0}$ for any $\mathbb{H}$-vector space $V$ in $Vect_\mathbb{H}$, and
$F(\psi\circ\phi)=F(\psi)\circ F(\phi)$ for any $\phi:V\rt W$ and $\psi:W\rt U$ in $Vect_\mathbb{H}$. The functor $F$ is clearly covariant and faithful. So $Vect_\mathbb{H}$ can be viewed as a subcategory of $Vect_\mathbb{C}$.

 Let $V$ be an $n$-dimensional right $\mathbb{H}$-vector space. Let $End_\mathbb{H}(V)$ be the set of $\mathbb{H}$-linear operators of $V$, which is an (associative) $\mathbb{R}$-algebra. One knows that the center of $End_\mathbb{H}(V)$ is $\{\lambda\cdot I|\lambda\in \mathbb{R}\}$, which isomorphic to $\mathbb{R}$.\\[3mm]

3. Now let us introduce the quaternionic structure on a complex vector space.

\begin{defi}
A conjugate linear map $\zeta$ from a complex vector space $W$ to itself satisfying $\zeta^2=-I$ is called a \textbf{quaternionic structure} on $W$.
\end{defi}

Let $V$ be a finite dimensional right $\mathbb{H}$-vector space. The map $r_j:V\rt V, v\mapsto v\cdot j$ gives rise to a map $J: V_0\rt V_0, v\mapsto v\cdot j$, which is conjugate linear and $J^2=-I$ because of (\ref{a}), thus $J$ is a quaternionic structure on $V_0$.

\bpr
If $\zeta$ is a quaternionic structure on a (finite-dimensional) complex vector space $W$, then \be \lb{e}W\times \mathbb{H}\rt W, (w,x+jy)\mapsto w.(x+jy)= wx+\zeta(w)y\ee (where $x,y\in \mathbb{C}$) defines an $\mathbb{H}$-vector space structure on $W$.
\epr
\bp
It is easy to verify that $v.1=v$, $(v+u).p=v.p+u.p$, and  \be\lb{m} w.(p+q)=w.p+w.q\ee for any $v,u,w\in W, p,q\in \mathbb{H}$.
For any $z\in \mathbb{C}$, let $r_z:W\rt W, w\mapsto w\cdot z$. One verifies $r_z \zeta=\zeta r_{\bar{z}}$ for any $z\in \mathbb{C}$ as follows: \bee (r_z \zeta)(w)= \zeta(w)z =\zeta(w\bar{z})=\zeta (r_{\bar{z}}(w))=(\zeta r_{\bar{z}})(w)\eee
Since \[r_z \zeta=\zeta r_{\bar{z}}~~ \forall z\in \mathbb{C} ~~\textrm{and}~~ \zeta^2=-I,\]
By (\ref{m}) and (\ref{a}), it can be directly verified that $ w.(pq)=(w.p).q$ for any $w\in W, p,q\in \mathbb{H}$.
So (\ref{e})  defines a right $\mathbb{H}$-module structure on $W$. As $\mathbb{H}$ is a division ring and  $W$ is finitely generated $\mathbb{H}$-module (as $W$ is a finite dimensional complex vector space), $W$ is a free $\mathbb{H}$-module, i.e., an $\mathbb{H}$-vector space.
\ep

A quaternionic  vector space can be viewed as a complex vector space  with a quaternionic structure. And,  a linear map between quaternionic  vector spaces  can be viewed as a linear map between complex vector spaces with quaternionic structures which commutes with the quaternionic structures.

\bpr\lb{j}
Assume that $U$ is a finite dimensional complex vector space with a quaternionic structure $\zeta$.
Then

(1) $U$ has no 1-dimensional $\zeta$-invariant subspace.

(2) $U=\bigoplus_{i=1}^k U_i$ where each $U_i$ is a 2-dimensional $\zeta$-invariant subspace. In particular, the dimension of $U$ is even.
\epr
\bp
(1) If $U$ has a 1-dimensional $\zeta$-invariant subspace spanned by a nonzero vector $u$, then $\zeta(u)=u\cdot\lambda$ for some $\lambda\in \mathbb{C}$. Since $\zeta$ is conjugate linear and $\zeta^2=-I$, \[-u=\zeta^2(u)=\zeta(\zeta u)=\zeta(u\cdot\lambda)=\zeta(u)\cdot \bar{\lambda}=u\cdot \lambda\bar{\lambda},\] forcing $\lambda\bar{\lambda}=-1$, a contradiction.

(2) By the above proposition, $U$ can be regarded as a finite-dimensional $\mathbb{H}$-vector space. Thus $U$ can be written as a direct sum of 1-dimensional $\mathbb{H}$-subspaces. Then as a $\mathbb{C}$-vector space, $U$ can be written as a direct sum of 2-dimensional $\zeta$-invariant subspaces.
 \ep

4.  Let $V$ be an $n$-dimensional right $\mathbb{H}$-vector space and $\phi\in End_\mathbb{H}(V)$. If $v\in V$ is nonzero and for some $\lambda\in \mathbb{H}$ one has
\[\phi(v)=v\cdot \lambda\] then we call $\lambda$ an \textit{eigenvalue} of $\phi$ and $v$ a corresponding \textit{eigenvector}. In this case, $\phi(v\cdot p)=(v\cdot p) p^{-1}\lambda p$ for any $p\in \mathbb{H}\setminus \{0\}$. Thus $v\cdot \mathbb{H}$, the subspace of $V$ spanned by $v$, correspond to $[\lambda]$, the set of quaternions similar to $\lambda$. In particular, there exists some nonzero $u\in v\cdot \mathbb{H}$ such that $\phi(u)=u\cdot \mu$ with $\mu\in \mathbb{C}$.
\ble
Let $\lambda\in \mathbb{C}$. If $v\in V$ is an eigenvector of $\phi$ corresponding to $\lambda$, then ${v\cdot j}\in V$ is an eigenvector of $\phi$ corresponding to $\bar{\lambda}$.
\ele
\bp
One has \[\phi(v\cdot j)=\phi(v)\cdot j=v\cdot \lambda j=v\cdot j\bar{\lambda}=(v\cdot j)\cdot\bar{\lambda}.\]
\ep
 A matrix $Q=[q_{ij}]\in M_n(\mathbb{H})$ can be view as a linear map $\mathbb{H}^n\rt \mathbb{H}^n, v\mapsto Qv$, and its eigenvalues are usually called right eigenvalues of $Q$. See \cite{z1} and etc. It is clear that similar matrices have the same eigenvalues.\\[3mm]

5. Let $V$ be an $n$-dimensional right $\mathbb{H}$-vector space with basis $\mathcal{B}=\{v_1,\cdots,v_n\}$. In this paper $\mathbb{H}^n$ will always be the  right $\mathbb{H}$-vector space of $n$-column vectors and $M_n(\mathbb{H})$ be the $\mathbb{R}$-algebra of $n\times n$ quaternionic matrices

 \ble
 The coordinate map $V\rt \mathbb{H}^n, v_1\cdot a_1+\cdots+v_n\cdot a_n\mapsto [a_1,\cdots,a_n]^T$ is an isomorphism of $\mathbb{H}$-vector spaces.
 \ele

  Let $\phi\in End_\mathbb{H}(V)$.

 Assume $\phi(v_j)=\sum_{i=1}^n v_i\cdot q_{ij}$, with $q_{ij}\in \mathbb{H}$. Then the matrix of $\phi$ in the (ordered) basis $\mathcal{B}$ is defined to be \[M(\phi,\mathcal{B})=[q_{ij}]\in M_n(\mathbb{H})\] as usual. Then one  has that

 \[End_\mathbb{H}(V)\rt M_n(\mathbb{H}), \phi\mapsto M(\phi, \mathcal{B})\] is an $\mathbb{R}$-algebra isomorphism.

  As \[\mathbb{H}=1\cdot \mathbb{C}\oplus j\cdot \mathbb{C}\] one has

\ble
The set $\{v_1,\cdots,v_n\}$ of vectors in $V$ is a  basis of $V$ if and only if $\{v_1,\cdots,v_n,v_1\cdot j,\cdots,v_n\cdot j\}$ is a  basis of $V_0$.
\ele

If $\mathcal{B}=\{v_1,\cdots,v_n\}$  is a  basis of $V$, then we will denote the corresponding basis $\{v_1,\cdots,v_n,v_1\cdot j,\cdots,v_n\cdot j\}$  of $V_0$ by $\mathcal{B}_0$.

Let $\mathcal{B}=\{e_1,\cdots,e_n\}$ be the standard basis of $\mathbb{H}^n$, then $\mathcal{B}_0=\{e_1,\cdots,e_n,e_1\cdot j,\cdots,e_n\cdot j\}$ is a basis of $(\mathbb{H}^n)_0$, the space $\mathbb{H}^n$ viewed as a $\mathbb{C}$-vector space. The coordinate of $[x_1+j y_1,\cdots, x_n+j y_n]^T$ (with $x_n,y_n\in \mathbb{C}$) in $\mathcal{B_0}$ is $[x_1,\cdots, x_n,y_1,\cdots, y_n]^T$. This induces a $\mathbb{C}$-linear isomorphism:
 \[\mathbb{H}^n\rt \mathbb{C}^{2n}, [x_1+j y_1,\cdots, x_n+j y_n]^T\mapsto [x_1,\cdots, x_n,y_1,\cdots, y_n]^T.\]

Under this linear isomorphism, the map $r_j:\mathbb{H}^n\rt \mathbb{H}^n, v\mapsto v\cdot j$ maps $[x_1+j y_1,\cdots, x_n+j y_n]^T$ to $ [-\bar{y_1}+j \bar{x_1},\cdots, -\bar{y_n}+j \bar{x_n}]^T$. It induces the conjugate linear map

 $$J:\mathbb{C}^{2n}\rt \mathbb{C}^{2n},$$
 $$[x_1,\cdots, x_n,y_1,\cdots, y_n]^T \mapsto [-\bar{y_1},\cdots, -\bar{y_n}, \bar{x_1},\cdots, \bar{x_n}]^T,$$
 which is the standard quaternionic structure on $\mathbb{C}^{2n}$. For $v\in \mathbb{C}^{2n}$, $J(v)$ is called the\textit{ adjoint vector} of $v$ in \cite{zw}. This map $J$ plays an important role in the proof of the Jordan canonical form theorem  for quaternionic matrices in \cite{zw}.

 Let $V, \mathcal{B}$ and $\phi$ be as above. If we view $V$ as a $\mathbb{C}$-vector space, $V_0$, then $\phi$ defines a $\mathbb{C}$-linear operator on $V_0$ denoted by $\phi_0:V_0\rt V_0$.
 \bpr
Assume $M(\phi, \mathcal{B})=Y+Z\cdot j$ with $Y=[y_{st}],Z=[z_{st}]$ in $M_n(\mathbb{C})$. Then the matrix of $\phi_0$ in the basis $\mathcal{B}_0$ is \begin{eqnarray*}
              M(\phi_0,\mathcal{B}_0)= \begin{bmatrix}
               Y&-Z\\
            \bar{Z}& \bar{ Y}
        \end{bmatrix}.
    \end{eqnarray*}
 \epr
\bp
One has $\phi(v_t)=\sum_{s=1}^n v_s(y_{st}+z_{st}j)$. Then for $t=1,\cdots, n$
\bee\begin{split}
\phi_0(v_t)=\phi(v_t)&=\sum_{s=1}^n v_s(y_{st}+z_{st}j)\\
&=\sum_{s=1}^n v_s y_{st}+\sum_{s=1}^n (v_s\cdot j) \overline{z_{st}}
 \end{split}\eee  and
$$\phi_0(v_t\cdot j)=\phi(v_t\cdot j)=\phi(v_t)\cdot j.$$
So \bee\begin{split}
\phi_0(v_t\cdot j)=\phi(v_t)\cdot j&=\sum_{s=1}^n v_s(y_{st}+z_{st}j)\cdot j\\
&=\sum_{s=1}^n v_s (-z_{st}+y_{st}j)\\
&=\sum_{s=1}^n v_s (-z_{st})+\sum_{s=1}^n (v_s\cdot j) \overline{y_{st}}.
 \end{split}\eee

 Thus  the matrix of $\phi_0$ in the basis $\mathcal{B}_0$ is \begin{eqnarray*}
              M(\phi_0,\mathcal{B}_0)= \begin{bmatrix}
               Y&-Z\\
            \bar{Z}& \bar{ Y}
        \end{bmatrix}.
    \end{eqnarray*}
\ep

Note that if the basis of $V_0$ corresponding to $\mathcal{B}$ is chosen to be $\mathcal{B}_0^{'}=\{v_1,\cdots,v_n,-v_1\cdot j,\cdots,-v_n\cdot j\}$, then \begin{eqnarray*}
              M(\phi_0,\mathcal{B}_0^{'})= \begin{bmatrix}
               Y&Z\\
            -\bar{Z}& \bar{ Y}
        \end{bmatrix},
    \end{eqnarray*} which is called the\textit{ complex adjoint matrix} of $Y+Z\cdot j$ in \cite{z1}, \cite{zw} and etc, and plays a central role in the study of
    quaternionic square matrices.
\bpr
Let $\mathcal{B}$ be a basis of $V$ and $\mathcal{B}_0$ the corresponding basis of $V_0$.
The maps $\chi:End_\mathbb{H}(V)\rt End_\mathbb{C}(V_0),\phi\mapsto \phi_0$ and
\[\widetilde{\chi}:M_n(\mathbb{H})\rt M_{2n}(\mathbb{C}), M(\phi, \mathcal{B})\mapsto M(\phi_0, \mathcal{B}_0)\] that maps $ Y+Zj$ to $ \begin{bmatrix}
               Y&-Z\\
            \bar{Z}& \bar{ Y}
        \end{bmatrix}$
are both injective $\mathbb{R}$-algebra homomorphisms. Furthermore, one has the following commutative diagrams:
\epr
$$
\begin{array}[c]{ccc}
End_\mathbb{H}(V)&\stackrel{\cong}{\longrightarrow}&M_n (\mathbb{H})\\
\downarrow\scriptstyle{\chi}&&\downarrow\scriptstyle{\widetilde{\chi}}\\
End_\mathbb{C}(V_0)&\stackrel{\cong}{\longrightarrow}&M_{2n} (\mathbb{C})
\end{array}
$$

\bp
It is clear that $\chi$ is injective and easy to verify that for $\phi,\psi\in End_\mathbb{H}(V)$ and $\la\in \mathbb{R}$, $\chi(\phi+\psi)=\chi(\phi)+\chi(\psi), \chi(\phi\psi)=\chi(\phi)\chi(\psi)$ and $\chi(\la\phi)=\la\chi(\phi)$. Thus $\chi$ is an injective $\mathbb{R}$-algebra homomorphism. The commutativity of the diagram is also obvious.

Since the two horizontal maps are clearly $\mathbb{R}$-linear isomorphisms, one knows that $\widetilde{\chi}$ is also an injective $\mathbb{R}$-algebra homomorphism.
\ep
Let $L$ be the image of $\chi$. Then it is clear that $$L=\{\phi\in End_\mathbb{C}(V_0)|\phi J=J\phi\},$$ and $\chi:End_\mathbb{H}(V)\rt L$ is  a bijection.

 A linear operator $\phi\in End_\mathbb{H}(V)$ is called \textit{semisimple}, or \textit{diagonalizable}, if there is some basis of $V$ in which the matrix of $\phi$ is a diagonal matrix. It is clear that $\phi$ is semisimple if and only if $\phi$ has $n$ linearly independent eigenvectors, where $n=dim~V$.  A linear operator $\phi\in End_\mathbb{H}(V)$ is called \textit{nilpotent} if there exists some positive integer  $k$, $\phi^k=0$.   A linear operator $\phi\in End_\mathbb{H}(V)$ is called \textit{unipotent} if $\phi-I$ is nilpotent. The notion of semisimple (resp. nilpotent, unipotent) operators are similarly defined for operators on $V_0$.

 \bpr\lb{k}
 The map $\chi: End_\mathbb{H}(V)\rt L$ sets up a 1-1 correspondence between semisimple (resp. nilpotent, unipotent) elements of $End_\mathbb{H}(V)$ and $L$.
\epr
\bp
Assume $\phi\in End_\mathbb{H}(V)$. Then $\phi_0=\chi(\phi)\in End_\mathbb{C}(V_0)$.

If $\phi$ is semisimple, then $\phi$ has $n$ linearly independent eigenvectors $v_1,\cdots,v_n$ with the respective eigenvalues $\la_1,\cdots,\la_n\in \mathbb{H}$. By the discussion above, one can assume that  $\la_1,\cdots,\la_n\in \mathbb{C}$. Then $v_1,\cdots,v_n, J(v_1),\cdots,J(v_n)$ are $2n$ linearly independent eigenvectors of $\phi_0$ thus $\phi_0$ is semisimple.

 Conversely, if $\phi_0$ is semisimple, then  $\phi_0$ has $2n$ linearly independent eigenvectors $u_1,\cdots,u_{2n}$. Thus $V_0=\oplus_{i=1}^{2n} u_i \mathbb{C}$ and $V=\sum_{i=1}^{2n} u_i \mathbb{H}$. One can find $n$ vectors in $u_1,\cdots,u_{2n}$, assuming them to be $u_1,\cdots,u_{n}$, such that $V=\oplus_{i=1}^{n} u_i \mathbb{H}$. So $\phi$ is semisimple.

It is obvious that $\phi$ is nilpotent if and only if $\phi_0$ is nilpotent.

As $\chi$ maps the identity operator on $V$ to the identity operator on $V_0$,  $\phi$ is unipotent if and only if $\phi_0$ is  unipotent.

\ep

7. Let $V$ be an $n$-dimensional right $\mathbb{H}$-vector space. Assume $\mathcal{B}=\{v_1,\cdots,v_n\}$ and $\mathcal{B}^{'}=\{v_1^{'},\cdots,v_n^{'}\}$ be two bases of $V$. Let $\phi\in End_\mathbb{H}(V)$. Let $M(\phi,\mathcal{B})$ and $M(\phi,\mathcal{B^{'}})$ be the respective matrices of $\phi$ in $\mathcal{B}$ and $\mathcal{B^{'}}$. Assume \[v_j^{'}=\sum_{i=1}^n v_i\cdot p_{ij}, j=1,\cdots, n.\] Then the \textit{transition matrix} from $\mathcal{B}$ to $\mathcal{B^{'}}$ is defined to be  $P=[p_{ij}]\in M_n(\mathbb{H})$.

Then one also has \be\lb{s} M(\phi,\mathcal{B^{'}})=P^{-1}M(\phi,\mathcal{B})P\ee as in the usual linear algebra, and the proof is the same.\\[3mm]

8. Let $V$ be an $n$-dimensional right $\mathbb{H}$-vector space. Let $\phi\in End_\mathbb{H}(V)$. Then the kernel of $\phi$, $Ker(\phi)$, and the image of $\phi$, $Im(\phi)$ or $\phi(V)$, are both $\mathbb{H}$-subspaces of $V$.

The rank of $\phi$, $rank(\phi)$, is defined to be the dimension of $Im(\phi)$. Then one has

\[dim~Ker(\phi)+rank(\phi)=dim~V.\]

\section{Exponential maps and characteristic polynomials}
 \setcounter{equation}{0}\setcounter{theorem}{0}
 Let $V$ be an $n$-dimensional right $\mathbb{H}$-vector space, and $V_0$ the corresponding $2n$-dimensional right $\mathbb{C}$-vector space with the quaternionic structure $J$. Let $ GL(V_0)$ (resp. $ GL(V)$) be the real Lie group of invertible linear operators in $End_\mathbb{C}(V_0)$ (resp. in $End_\mathbb{H}(V)$). Identify $End_\mathbb{C}(V_0)$ (resp. $End_\mathbb{H}(V)$) with the Lie algebra of $GL(V_0)$ (resp. $ GL(V)$).

 Recall $$L=\{\phi\in End_\mathbb{C}(V_0)|\phi J=J\phi\},$$ which is a subalgebra of the real associative algebra $End_\mathbb{C}(V_0)$, and is also a real Lie algebra under the usual Lie bracket. Let $$G=\{\phi\in GL(V_0)|\phi J=J\phi\},$$ which is a closed Lie subgroup of $GL(V_0)$.

 Let \[exp_1: M_n(\mathbb{H})\rt M_n(\mathbb{H}), A\mapsto \sum_{k=0}^{\infty} \frac{1}{k!}A^k\] be the exponential map on $M_n(\mathbb{H})$. It can be directly verified that the series $\sum_{k=0}^{\infty} \frac{1}{k!}A^k$ is convergent \cite{t} and $exp_1 (A)\in GL_n(\mathbb{H})$.

  Similarly, let \[exp_2: End_\mathbb{H}(V)\rt GL(V), \phi\mapsto \sum_{k=0}^{\infty} \frac{1}{k!}\phi^k\] be the exponential map on $End_\mathbb{H}(V)$. Fix some basis $\mathcal{B}$ of $V$, then one has the following commutative diagram:
  $$
\begin{array}[c]{ccc}
GL(V)&\longrightarrow& GL_n(\mathbb{\mathbb{H}})\\
\uparrow\scriptstyle{exp_1}   &  &  \uparrow\scriptstyle{exp_2}        \\
End(V) &\longrightarrow&M_n(\mathbb{\mathbb{H}}) \\
\end{array}
$$
\begin{center}
$\text{Figure 1}$
\end{center}

Here the horizontal maps take $\phi\in GL(V)$ (resp. $\phi\in End(V)$) to $M(\phi,\mathcal{B})$.

For any $\psi\in End_\mathbb{C}(V_0)$, define $\xi(\psi)=J\psi J^{-1}$, which equals $-J\psi J$. Then $\xi(\psi)$ is also in $End_\mathbb{C}(V_0)$.  Let $$End_\mathbb{C}(V_0)^\xi=\{\psi\in End_\mathbb{C}(V_0)|\xi(\psi)=\psi\}$$ be the set of operators fixed by $\xi$. It is clear that $$End_\mathbb{C}(V_0)^\xi=L.$$

\ble  The map $\xi:End_\mathbb{C}(V_0)\rt End_\mathbb{C}(V_0)$ is  conjugate linear.
One has $\xi^2=Id$, $\xi(\la I)=\bar{\la}I$ for any $\la\in \mathbb{C}$,  $\xi(\psi_1+\psi_2)=\xi(\psi_1)+\xi(\psi_2)$ and $\xi(\psi_1\psi_2)=\xi(\psi_1) \xi(\psi_2)$ for any $\psi_1,\psi_2\in End_\mathbb{C}(V_0)$.
\ele
\bp 
The fact that $\xi^2=Id$, $\xi(\psi_1+\psi_2)=\xi(\psi_1)+\xi(\psi_2)$ and $\xi(\psi_1\psi_2)=\xi(\psi_1) \xi(\psi_2)$ are directly verified.

As $J(\la I)J^{-1}=(\bar{\la}I) J J^{-1}=\bar{\la}I $, one has $\xi(\la I)=\bar{\la}I$.

Then for any $\psi\in End_\mathbb{C}(V_0)$, $\xi(\la\psi)=\xi(\la I\cdot\psi)=\xi(\la I)\xi(\psi)=(\bar{\la} I)\xi(\psi)=\bar{\la} \xi(\psi)$. Thus $\xi$ is a conjugate linear operator on $ End_\mathbb{C}(V_0)$.
\ep
As $\xi^2=Id$ and $\xi$ is conjugate linear, $L=End_\mathbb{C}(V_0)^\xi$ is a real subspace of $End_\mathbb{C}(V_0)$ and $End_\mathbb{C}(V_0)=L\oplus i\cdot L$.

Let $\phi\in End_\mathbb{H}(V)$. Then $\phi_0\in L$. Let $\Gamma(\phi_0)$ be the set of eigenvalues of $\phi_0$ in $\mathbb{C}$.

 If $\la$ is an an eigenvalue of $\phi_0$, then \[V_0(\la)=Ker~(\phi_0-\la I)^{2n}\] is called the\textit{ generalized eigenspace} of $\phi_0$ belonging to $\la$. Let $m(\la)=dim~V_0(\la)$, which is the multiplicity of the eigenvalue $\la$.

 \ble
Assume that $\la\in \Gamma(\phi_0)$. Then

(1)$\bar{\la}\in \Gamma(\phi_0)$, and $J$ maps  $V_0(\la)$ isomorphically onto $V_0(\bar{\la})$;

(2) $m(\la)=m(\bar{\la})$;

(3) If $\la\in \mathbb{R}$, then $m(\la)$ is even.
\ele
\bp

(1) Assume that $\phi_0(v)=v\cdot \la$ with $v$ a nonzero vector in $V_0$. Then $\phi_0(J v)=J \phi_0(v)=J(v\la)=Jv\cdot \bar{\la}$ with $Jv$ a nonzero vector. Thus $\bar{\la}\in \Gamma(\phi_0)$.

 For any $\phi_0\in L$ and $\la\in \Gamma(\phi_0)$, $\xi [(\phi_0-\la I)^{2n}]=[\xi (\phi_0-\la I)]^{2n}=(\phi_0-\bar{\la} I)^{2n}$. So $J$ maps  $V_0(\la)=Ker (\phi_0-\la I)^{2n}$ isomorphically onto $V_0(\bar{\la})=Ker (\phi_0-\bar{\la} I)^{2n}$.

(2) follows from (1).

(3)  If $\la\in \mathbb{R}$ then it follows from (1) that $V_0(\la)$ is $J$-invariant. By Proposition \ref{j}, $m(\la)=dim~V_0(\la)$ is even.
\ep

By this lemma, one can assume that $\Gamma(\phi_0)=\{\la_1,\bar{\la_1},\cdots,\la_k,\bar{\la_k}\}\cup \{\mu_1,\cdots,\mu_t\}$, where $\mu_i\in \mathbb{R}$ and $\la_i\in \mathbb{C}\setminus \mathbb{R}$,  $m(\la_i)=m(\bar{\la_i})$ for any $i$ and  $m(\mu_j)$ is even for any $j$.

For any $\phi\in End_\mathbb{H}(V)$, motivated by \cite{z1}, we define its \textit{characteristic polynomial} to be $p(x)=det(xI-\phi_0)$.

\ble
One has $p(x)\in \mathbb{R}[x]$.
\ele
\bp
By the above lemma, one has $$p(x)=\prod_{i=1}^k[(x-\la_i)(x-\bar{\la_i})]^{m(\la_i)}\prod_{j=1}^t (x-\mu_j)^{m(\mu_j)}\in \mathbb{R}[x].$$
\ep
 If the trace  of $\phi$ and the determinant of $\phi$ are defined to be the trace and determinant of $\phi_0$ respectively, then it is easy to see that they are both real.

Here is the Cayley-Hamilton theorem for quaternionic linear operators.
\bthm  \cite{z1}
Let $p(x)$ be the characteristic polynomial of $\phi\in End_\mathbb{H}(V)$. Then  $p(\phi)=0$.
\ethm
\bp
Recall the injective $\mathbb{R}$-algebra homomorphism $\chi:End_\mathbb{H}(V)\rt End_\mathbb{C}(V_0)$. As $p(x)\in \mathbb{R}[x]$,  $\chi(p(\phi))=p(\chi(\phi))=p(\phi_0)=0$. Since $\chi$ is injective, one has $p(\phi)=0$.
\ep
\bigskip
For any $\psi\in GL(V_0)\subseteq End_\mathbb{C}(V_0)$,  $\xi(\psi)=J\psi J^{-1}$ is also in $GL(V_0)$.  The restriction of $\xi$ to $GL(V_0)$,  denoted by $\widetilde{\xi}:GL(V_0)\rt GL(V_0)$, is clearly a Lie group homomorphism. Then the differential of the Lie group homomorphism $\widetilde{\xi}:GL(V_0)\rt GL(V_0)$ is just  $\xi:End_\mathbb{C}(V_0)\rt End_\mathbb{C}(V_0)$. As $ GL(V_0)^{\widetilde{\xi}}=G$ and $ End_\mathbb{C}(V_0)^\xi=L$. One gets

 \ble\lb{h}
 The Lie algebra of $G$ is $L$. And, $exp:End_\mathbb{C}(V_0)\rt GL(V_0)$ maps $L$ into $G$.
 \ele


The images of $End_\mathbb{H}(V)$ and $GL(V)$ under the map $\chi:End_\mathbb{H}(V)\rt End_\mathbb{C}(V_0)$ are respectively $\chi(End_\mathbb{H}(V))=L$ and $ \chi(GL(V))=G$. So one has the following commuting diagram\\[2mm]

 $$
\begin{array}[c]{ccccc}
GL(V)&\xrightarrow[\cong]{\chi }& G &\xrightarrow{} &GL(V_0)\\
\uparrow\scriptstyle{exp}   &  & \quad \uparrow\scriptstyle{exp}        &   &\uparrow\scriptstyle{exp}\\
End_\mathbb{H}(V) &\xrightarrow[\cong]{\chi}&L &\xrightarrow{}& End_{\mathbb{C}}(V_0)
\end{array}
$$

\begin{center}
$\text{Figure 2}$
\end{center}

\section{Proof of Jordan-Chevalley decomposition for quaternionic linear operators}
 \setcounter{equation}{0}\setcounter{theorem}{0}
 Assume that $V$ is an $n$-dimensional right $\mathbb{H}$-vector space and $\phi\in End_\mathbb{H}(V)$. Let $V_0$ be the corresponding $2n$-dimensional right $\mathbb{C}$-vector space and $\phi_0\in End_\mathbb{C}(V_0)$ the corresponding $\mathbb{C}$-linear operator induced by $\phi$.

1. The additive Jordan-Chevalley decomposition of operators in $End_\mathbb{H}(V)$ and in $L$

\bthm\lb{f}
Let $\psi\in L$. Then

(1) $\psi$ can be uniquely written as $\psi=\psi_s+\psi_n$ with $\psi_s,\psi_n\in L$, where $\psi_s$ is semisimple, $\psi_n$ is nilpotent, and $\psi_s\psi_n=\psi_n\psi_s$.

(2) There exist polynomials $f(x), g(x)\in \mathbb{R}[x]$ without constant term such that $\psi_s=f(\psi)$ and $\psi_n=g(\psi).$
\ethm
\bp
Denote $V_0$ by $U$. Assume that the set of eigenvalues of $\phi_0$ is $\Gamma(\phi_0)=\{\la_1,\bar{\la_1},\cdots,\la_k,\bar{\la_k}\}\cup \{\mu_1,\cdots,\mu_t\}$, where $\mu_i\in \mathbb{R}$ and $\la_i\in \mathbb{C}\setminus \mathbb{R}$. One has \bee U=\oplus_{i=1}^k (U(\la_i)\oplus U(\bar{\la_i}))\bigoplus (\oplus_{j=1}^t U(\mu_j))\eee as the decomposition of $U$ into generalized eigenspaces of $\psi$.

One has $det (x\cdot I-\psi)=\Pi_{i=1}^k(x-\lambda_i)^{m(\la_i)}(x-\bar{\lambda_i})^{m(\la_i)}\cdot\Pi_{j=1}^t(x-\mu_j)^{m(\mu_j)}$,  where $\mu_i\in \mathbb{R}$ and $\la_i\in \mathbb{C}\setminus \mathbb{R}$. Consider the following system of congruences for $h(x)\in \mathbb{C}[x]$:
$$
\left\{
\begin{aligned}
h(x)\equiv\lambda_i  ~~~(mod~&(x-\lambda_i)^{m(\la_i)}), ~~~i=1,\cdots,k \\
h(x)\equiv\bar{\lambda_i}  ~~~(mod~&(x - \bar{\lambda_i})^{m(\la_i)}),~~~ i=1,\cdots,k\\
h(x)\equiv \mu_j~~~ (mod~&(x-\mu_j)^{m(\mu_j)}), ~~~j=1,\cdots,t\\
h(x)\equiv 0 ~~~ (mod~ & x).
\end{aligned}
\right.
$$
Notice that the last congruence is superfluous if 0 is an eigenvalue of $\psi$. Now apply the Chinese Remainder Theorem for $\mathbb{C}[x]$. Since the congruences in the system have pairwise relatively prime moduli, there is some solution $h_0(x)$ for it.

Notice that if $h_0(x)\equiv\lambda_i~~(mod(x-\lambda_i)^{n_i})$, then $h_0(x)=(x-\lambda_i)^{n_i}\cdot k(x)+\lambda_i$ and
 $\overline{h_0(x)}=(x-\bar{\lambda_i})^{n_i}\cdot\overline{k(x)}+\bar{\lambda_i}$ for some $k(x)\in \mathbb{C}[x]$, so
 $\overline{h_0(x)}\equiv\bar{\lambda_i}~~(mod(x-\bar{\lambda_i})^{n_i})$\\
Similarly,we can easily prove that $\overline{h_0(x)}$ satisfies each of the congruences. Let $f(x)=\frac{1}{2}(h_0(x)+\overline{h_0(x)})$,it is easily seen that $f(x)\in \mathbb{R}[x]$ is also the solution of the congruence equations. Then,
$$
f(\psi)|_W=
\left\{
\begin{aligned}
\lambda_i \cdot I ,~~~W=U(\lambda_i );\\
\bar{\lambda_i} \cdot I ,~~~W=U(\bar{\lambda_i});\\
\mu_j \cdot I,~~~W=U(\mu_j).
\end{aligned}
\right.
$$
Let $\psi_s= f(\psi)$. Then $\psi_s$ is semisimple and is in $L$. let $g(x)=x-f(x)$ and $\psi_n =g(\psi)$. It is clear that  $\psi_n$
is nilpotent and is in $L$. One has $\psi=\psi_s+\psi_n$ and $\psi_s\psi_n=\psi_n\psi_s$. Obviously,$f(x)$ and $g(x)$ have no constant term.

It remains only to prove the uniqueness assertion in (1). Assume $\psi=\psi^{'}_s+\psi^{'}_n$ is another such decomposition, where $\psi^{'}_s, \psi^{'}_n\in L$, $\psi^{'}_s$ is semisimple, $\psi^{'}_n$ is nilpotent, and $\psi^{'}_s\psi^{'}_n=\psi^{'}_n\psi^{'}_s$. It is clear that $\psi^{'}_s$ and $\psi^{'}_n$ commute with $\psi$.

Since $\psi_s=f(\psi),\psi_n=g(\psi)$ ,$\psi^{'}_s$ and $\psi^{'}_n$ commute with $\psi_s$ and $\psi_n$. As $\psi_s+\psi_n=\psi^{'}_s+\psi^{'}_n$, $\psi_s-\psi^{'}_s=-\psi_n+\psi^{'}_n$, which are  both semisimple and nilpotent, thus must be the zero operator. This forces that $\psi_s=\psi^{'}_s,\psi_n=\psi^{'}_n$.
\ep

Now we can show the additive Jordan-Chevalley decomposition for quaternionic linear operators.
\bco \lb{g}
Let $\phi\in End_\mathbb{H}(V)$. Then

(1) $\phi$ can be uniquely written as $\phi=\phi_s+\phi_n$ with $\phi_s,\phi_n\in End_\mathbb{H}(V)$, where $\phi_s$ is semisimple, $\phi_n$ is nilpotent, and $\phi_s\phi_n=\phi_n\phi_s$.

(2) There exist polynomials $f(x), g(x)\in \mathbb{R}[x]$ without constant term such that $\phi_s=f(\phi), \phi_n=g(\phi).$
\eco
\bp
Let $\phi_0\in L$ be the operator on $V_0$ corresponding to $\phi$ and apply last theorem to $\phi_0$. Then there exists $(\phi_0)_s,(\phi_0)_n$ in $L$, $\phi_0=(\phi_0)_s+(\phi_0)_n$, $(\phi_0)_s$ is semisimple, $(\phi_0)_n$ is nilpotent, and $(\phi_0)_s(\phi_0)_n=(\phi_0)_n(\phi_0)_s$. And, there exist polynomials $f(x), g(x)\in \mathbb{R}[x]$ without constant term such that $(\phi_0)_s=f(\phi_0), (\phi_0)_n=g(\phi_0).$

Consider the $\mathbb{R}$-algebra homomorphism $\chi:End_\mathbb{H}(V)\rt L $. One has $\chi(f(\phi))=f(\phi_0)=(\phi_0)_s$ is semisimple, so $f(\phi)$ is semisimple by Proposition \ref{k}.  One has $\chi(g(\phi))=g(\phi_0)=(\phi_0)_n$ is nilpotent, so $g(\phi)$ is nilpotent by Proposition \ref{k}. Then  $$\chi(f(\phi)+g(\phi))=(\phi_0)_s+ (\phi_0)_n=\chi(\phi).$$ As  $\chi$ is injective, $\phi=f(\phi)+g(\phi)$. Thus one has the desired  decomposition $\phi=\phi_s+\phi_n$  with  $\phi_s=f(\phi)$ and $\phi_n=g(\phi)$. As such decomposition of $\phi_0$ is unique and $\chi$ is injective, such decomposition for $\phi$ is also unique.

\ep

2. The multiplicative Jordan-Chevalley decomposition of operators in $GL(V)$ and in $G$


\bthm
Let $\psi\in G$. Then

(1) $\psi$ can be uniquely written as $\psi=\psi_s\psi_u$ with $\psi_s,\psi_u\in G$, where $\psi_s$ is semisimple, $\psi_u$ is unipotent, and $\psi_s\psi_u=\psi_u\psi_s$.

(2) There exist polynomials $f(x), h(x)\in \mathbb{R}[x]$ with $f(0)=0, h(0)=1$ such that $\psi_s=f(\psi), \psi_u=h(\psi).$
\ethm
\bp
As $G\subseteq L$, $\psi\in L$. By Theorem \ref{f}, $\psi$ can be uniquely written as $\psi=\psi_s+\psi_n$, $\psi_s,\psi_n\in L$, where $\psi_s$ is semisimple, $\psi_n$ is nilpotent, and $\psi_s\psi_n=\psi_n\psi_s$. And, there exist polynomials $f(x), g(x)\in \mathbb{R}[x]$ without constant term such that $\psi_s=f(\psi)$ and $\psi_n=g(\psi).$

As $\psi$ is invertible, $\psi_s$ is invertible. Let $\psi_u=I+\psi_s^{-1}\psi_n$, which is unipotent. Then $\psi=\psi_s+\psi_n=\psi_s(I+\psi_s^{-1}\psi_n)=\psi_s\psi_u$. One has
$\psi_s\psi_u=\psi_u\psi_s$. Let $k(x)=det(xI-\psi)\in \mathbb{R}[x]$ be the characteristic polynomial of $\psi$, which is just the characteristic polynomial of $\psi_s$ as $\psi$ and $\psi_s$ have the same set of eigenvalues. One has $k(0)\neq 0$ as $\psi$ is invertible, so there exists $q(x)\in \mathbb{R}[x],\psi_s^{-1}=q(\psi_s)$.

Let $h(x)=1+(q\circ f)(x)\cdot g(x)\in \mathbb{R}[x]$. Since $g(0)=0$, $h(0)=1$. Then $h(\psi)=I+q( f(\psi))\cdot g(\psi)=I+q(\psi_s)\cdot\psi_n=I+\psi_s^{-1}\psi_n=\psi_u$. The proof of (2) is complete. What's left to prove is the uniqueness assertion in (1).

 Assume that $\psi=\psi^{'}_s\psi^{'}_u$ is another such decomposition with $\psi^{'}_s$ semisimple, $\psi^{'}_u$ unipotent, and $\psi^{'}_s\psi^{'}_u=\psi^{'}_u\psi^{'}_s$. Because of (2), $\psi^{'}_s$ and $\psi^{'}_u$ commute with $\psi_s$ and $\psi_u$. Then $(\psi^{'}_s)^{-1}\psi_s=\psi^{'}_u\psi_u^{-1}$, which are both semisimple and unipotent thus must be $I$. This forces that $\psi_s=\psi^{'}_s,\psi_u=\psi^{'}_u$.
\ep

Below is the multiplicative Jordan-Chevalley decomposition for quaternionic linear operators. Its proof is similar to that of Corollary \ref{g}, and we omit it.
\bco Let $\phi\in GL(V)$. Then

(1) $\phi$ can be uniquely written as $\phi=\phi_s\phi_u$ with $\phi_s,\phi_u\in GL(V)$, where $\phi_s$ is semisimple, $\phi_u$ is unipotent, and $\phi_s\phi_u=\phi_u\phi_s$.

(2) There exist polynomials $f(x), h(x)\in \mathbb{R}[x]$ with $f(0)=0, h(0)=1$ such that $\phi_s=f(\phi), \phi_u=h(\phi).$
\eco

 The additive Jordan-Chevalley decomposition and multiplicative Jordan-Chevalley decomposition are related by the exponential map as follows.
 \bpr
For any $\phi\in End_\mathbb{H}(V)$, let $\phi=\phi_s+\phi_n$ be the additive Jordan-Chevalley decomposition, with $\phi_s$ semisimple and $\phi_n$ nilpotent and $\phi_s\phi_n=\phi_n\phi_s$. Let $\psi=exp(\phi)$, $\psi_s=exp(\phi_s)$ and $\psi_u=exp(\phi_n)$.  Then $\psi$ is in $GL(V)$, and $\psi= \psi_s \psi_u$ is the multiplicative Jordan-Chevalley decomposition for $\psi$, with $\psi_s$ semisimple and $\psi_u$ unipotent and $\psi_s\psi_u=\psi_u\psi_s$.
\epr
\bp
As $\phi_s\phi_n=\phi_n\phi_s$, $$\psi=exp(\phi)=exp(\phi_s+\phi_n)=exp(\phi_s)exp(\phi_n)=\psi_s \psi_u.$$ It is clear that $exp(\phi_s)$ is  semisimple, $exp(\phi_n)$ is unipotent, and $exp(\phi_s)exp(\phi_n)=exp(\phi_n)exp(\phi_s)$. Then it follows from the uniqueness of the  multiplicative Jordan-Chevalley decomposition that $\psi= \psi_s \psi_u$ is just the multiplicative Jordan-Chevalley decomposition for $\psi$.
\ep

It is clear that analogous result also holds for complex linear operators.
\section{Proof of the theorem of Jordan canonical form for quaternionic linear operators}
 \setcounter{equation}{0}\setcounter{theorem}{0}
  Recall that
 a $k\times k$ matrix of the following type
\begin{eqnarray*}
\begin{array}{c@{\hspace{-5pt}}l}
 \text{   } \mathrm{J_k(\la)}=
 \begin{bmatrix}
 \lambda&1&&\\
 &\ddots& \ddots&\\
 &&\ddots& 1\\
 &&&\lambda
 \end{bmatrix}
  \begin{array}{l}
          \rule{0mm}{5mm} \\
        \rule{0mm}{5mm} \\
        \rule{0mm}{5mm} \\
        \rule{0mm}{5mm} \\
        \end{array}
  \end{array}
\end{eqnarray*} is called a \textit{Jordan block}, where $\la\in \mathbb{H}$.

Assume that $U$ is a complex vector space and $\psi\in End(U)$. Let $S$ be the set of $\psi$-invariant subspace $W$ of $U$ such that with respect to some basis $\mathcal{B}$ of $W$, $M(\psi|_W, \mathcal{B})$ is a Jordan block. Then $S$ ordered by inclusion, i.e. $(S,\subseteq)$, is a partially ordered set.  We call maximal elements in $S$ \textbf{Jordan subspaces} of $\psi$. A basis $\mathcal{B}$ of a Jordan subspace with respect to which $M(\psi|_W, \mathcal{B})$ is a Jordan block is usually called a \textit{Jordan chain}. The Jordan canonical form theorem says that $U$ can be decomposed as a direct sum of the Jordan subspaces of $\psi$. Such decomposition is not unique (The authors mistakenly assumed that the decomposition is unique in \cite{jz}), but the corresponding set of Jordan blocks are unique up to permutation.

The following result is simple and its proof is omitted.
\ble\lb{o}
Assume that $W$ is a $\psi$-invariant subspace of $U$ such that with respect to some basis $\mathcal{B}=\{v_1,\cdots,v_k\}$ of $W$, $M(\psi|_W, \mathcal{B})$ is a Jordan block. Then the only nonzero $\psi$-invariant subspaces of $W$ are $W_i=span\{v_1,\cdots,v_i\}$, for $i=1,\cdots, k$.
\ele
Next assume that $U$ is a complex vector space with the quaternionic structure $\zeta$,  and $\psi$ is a linear operator on $U$ commuting with $\zeta$.

\bpr
Assume that $W\subseteq U$ is a  Jordan subspace of $\psi$. Then $\zeta W$ is also a Jordan subspace of $\psi$ and
$\zeta W\cap W=0$.

\epr
\bp Assume $M(\psi|_W, \mathcal{B})=\mathrm{J_k(\la)}$ for some basis $\mathcal{B}=\{v_1,\cdots,v_k\}$ of $W$. It is clear that $\zeta  W$ is also a Jordan subspace of $\psi$, with basis $\zeta \mathcal{B}$, and $M(\psi|_{\zeta W}, \zeta \mathcal{B})=\mathrm{J_k(\bar{\la})}$. If $\la\notin \mathbb{R}$, it is clear that $\zeta W\cap W=0$. Next assume that $\la\in \mathbb{R}$.

Assume $\zeta W\cap W\neq 0$. Then $\psi(\zeta W\cap W)\subseteq \zeta W\cap W$. So $\zeta W\cap W$ is a $\psi$-invariant subspace in $W$. By Lemma \ref{o}, $\zeta W\cap W=W_i$ for some $i=1,\cdots,k$, where $W_i=span\{v_1,\cdots,v_i\}$.

Since $(\psi-\la I)\zeta=\zeta(\psi-\la I)$, $W_{i-1}=(\psi-\la I)(W_i)$ is also $\zeta$-invariant. Similarly one has that  $W_{1}=(\psi-\la I)^{i-1}(W_i)$ is also $\zeta$-invariant. But $W_1$ is 1-dimensional and $\zeta$ has no 1-dimensional invariant subspace, which is a contradiction.
\ep
\bpr\lb{q}
 Assume that $\psi$ is nilpotent on $U$. Then there exist Jordan subspaces $W_i$ and $\zeta W_i$ of $\psi$, $i=1,\cdots,k$, such that $U=\bigoplus_{i=1}^k(W_i\oplus \zeta W_i)$.
\epr
\bp
We prove it by induction on $dim~U$.

Assume $dim~U=2m$. If $m=1$, then there is some nonzero $v\in U$ with $\psi(v)=0$. Then $\psi(\zeta v)=0$. One knows that $v$ and $\zeta v$ are linearly independent. Then  $U$ can be decomposed as the direct sum of the pair of 1-dimensional Jordan subspaces, spanned by $v$ and $\zeta v$ respectively.

Assume that $m>1$ and the results holds for spaces with dimension $<2m$. Since $\psi \zeta=\zeta\psi$, $U_1=\psi(U)$ is also $\zeta$-invariant and $\psi$-invariant. As $\psi$ is nilpotent on $U$, $dim~U_1<dim~U$. By induction, $U_1$ can be decomposed as a direct sum of $\zeta$-paired Jordan subspaces of $\psi|_{U_1}$, and denote them by $W_1,\cdots,W_k, \zeta W_1,\cdots,\zeta W_k$.

Assume that the Jordan chain of $W_i$ is $\{\psi^t (v_i)|t=l_i,l_i-1,\cdots,1,0\}$, where $\psi^{l_i+1} (v_i)=0$. As $v_i$ is in the image of $\psi$, there exists some $w_i\in U$ with $\psi(w_i)=v_i$. Then $\psi(\zeta w_i)=\zeta v_i$. Let $\widetilde{W_i}=W_i\oplus \langle w_i\rangle$, then $\zeta\widetilde{W_i}=\zeta W_i\oplus \langle \zeta w_i\rangle$. It is clear that $\widetilde{W_i}$ and $\zeta\widetilde{W_i}$ are Jordan subspaces of $\psi$ thus $\widetilde{W_i}\cap \zeta\widetilde{W_i}=0$. We will show that the union of the $2k$ Jordan chains
$$ \psi^{l_1+1}(w_1),\cdots,\psi (w_1),w_1,\cdots, \psi^{l_k+1}(w_k),\cdots,\psi( w_k),w_k$$
 $$\psi^{l_1+1}(\zeta w_1),\cdots,\psi (\zeta w_1),\zeta w_1,\cdots, \psi^{l_k+1}(\zeta w_k),\cdots,\psi(\zeta w_k),\zeta w_k$$~ are~ linearly~ independent.
Note that $\psi^{l_i+2}(w_i)=\psi^{l_i+2}(\zeta w_i)=0$. Assume that \be\lb{b}\sum_{i=1}^k(\sum_{s=0}^{l_i+1} \psi^s( w_i)\cdot a_{i,s})+\sum_{i=1}^k(\sum_{s=0}^{l_i+1} \psi^s(\zeta w_i)\cdot a_{i,s}^{'})=0\ee with $a_{i,s}, a_{i,s}^{'}\in \mathbb{C}$. Applying $\psi$ on both sides one gets that $$a_{i,0}=a_{i,1}=\cdots=a_{i,l_i}=0, a_{i,0}^{'}=a_{i,1}^{'}=\cdots=a_{i,l_i}^{'}=0$$ for $i=1,\cdots,k$.  Then (\ref{b}) becomes \bee\sum_{i=1}^k \psi^{l_i+1} (w_i)\cdot a_{i,l_i+1}+\sum_{i=1}^k \psi^{l_i+1} (\zeta w_i)\cdot a_{i,l_i+1}^{'}=\sum_{i=1}^k \psi^{l_i}( v_i)\cdot a_{i,l_i+1}+\psi^{l_i}(\zeta v_i)\cdot a_{i,l_i+1}^{'}=0,\eee which forces $a_{i,l_{i+1}}=a_{i,l_{i+1}}^{'}=0$ for $i=1,\cdots,k$. So we have proved that the union of the $2k$ Jordan chains are linearly independent. It follows that the sum of the subspaces $\widetilde{W_1},\cdots,\widetilde{W_k}, \zeta\widetilde{W_1},\cdots,\zeta\widetilde{W_k}$, call it $Y$, is a direct sum. Note that $\psi(Y)=\psi(U)$.
If $Y=U$ then it is done.

 Assume that $Y\subsetneqq U$. Choose some $v\in U\setminus Y$. Then there is some $w\in Y$ with $\psi (v)=\psi (w)$. Let $u_1=v-w$. Then $\psi (u_1)=0$ and $u_1\notin  Y$. Then $u_1,\zeta u_1$ are linearly independent. Let $\langle u_1,\zeta u_1\rangle$ denote $span\{u_1,\zeta u_1\}$. Then $Y+\langle u_1,\zeta u_1\rangle$ must be a direct sum as $Y\cap \langle u_1,\zeta u_1\rangle$ is $\zeta$-invariant and has dimension less than 2. If $Y\oplus \langle u_1,\zeta u_1\rangle$ is still properly contained in $U$, one can continue this way and finally write $U=Y\oplus \langle u_1,\zeta u_1\rangle\oplus \langle u_2,\zeta u_2\rangle\oplus \cdots\oplus \langle u_t,\zeta u_t\rangle$ since $U$ is finite dimensional. It is clear that each $\langle u_i\rangle$ (resp. $\zeta\langle u_i\rangle$) is a Jordan subspace of $\psi$ in $U$. Thus \bee U=(\widetilde{W_1}\oplus \zeta\widetilde{W_1})\oplus\cdots\oplus(\widetilde{W_k}\oplus \zeta\widetilde{W_k})\oplus(\langle u_1\rangle\oplus \zeta\langle u_1\rangle)\oplus\cdots\oplus(\langle u_t\rangle\oplus \zeta\langle u_t\rangle)\eee is a decomposition into a direct sum of $\zeta$-paired Jordan subspaces.
\ep

 Now we come back to the usual setting. Assume that $V$ is an $n$-dimensional right $\mathbb{H}$-vector space and $\phi\in End_\mathbb{H}(V)$. Let $V_0$ be the corresponding $2n$-dimensional right $\mathbb{C}$-vector space with the corresponding quaternionic structure $J$, and $\phi_0\in End_\mathbb{C}(V_0)$ the corresponding $\mathbb{C}$-linear operator induced by $\phi$.

We will show that
$V_0$ can be written as a direct sum as follows
\[V_0=\oplus_{i=1}^l (W_i\oplus J(W_i)),\] where $W_i$ and $J(W_i)$ are both  Jordan subspaces of $\phi_0$ for each $i$.

\bthm\lb{i}
The complex space $V_0$ can be decomposed as a direct sum of $J$-paired Jordan subspaces of $\phi_0$.
\ethm
\bp
Denote $V_0$ by $U$. Then $dim~U=2n$. Assume that the set  of eigenvalues of $\phi_0$ is $\Gamma(\phi_0)=\{\la_1,\bar{\la_1},\cdots,\la_k,\bar{\la_k}\}\cup \{\mu_1,\cdots,\mu_t\}$, where $\mu_i\in \mathbb{R}$ and $\la_i\in \mathbb{C}\setminus \mathbb{R}$.  Let $U(\la_i)$ and $U(\mu_j)$ be the corresponding generalized eigenspaces of $\phi_0$.  Then
\[U=\oplus_{i=1}^k(U(\la_i)\oplus U(\bar{\la_i}))\bigoplus(\oplus_{j=1}^t U(\mu_j)).\]

The space $U(\la_i)$ can be decomposed as a direct sum of the Jordan subspaces $W_{i,s}$ of $\phi_0$. Then $U(\bar{\la_i})=J(U(\la_i))$ is the direct sum of the Jordan subspaces $J (W_{i,s})$ of $\phi_0$ .

Now consider $U(\mu_j)$. Denote $\phi_0$ by $\psi$, let $\psi_j=\psi|_{U(\mu_j)}$, then $\psi_j-\mu_j I$ is nilpotent on $U(\mu_j)$. By Proposition  \ref{q}, $U(\mu_j)$ can also be decomposed as a direct sum of $J$-paired Jordan subspaces, $Q_s\oplus J (Q_s)$, of $\psi_j-\mu_j I$, with Jordan blocks $\mathrm{J_m(0)}\oplus \mathrm{J_m(0)}$. But such Jordan subspaces are also Jordan subspaces of $\psi_j$, with Jordan blocks $\mathrm{J_m(\la)}\oplus \mathrm{J_m(\la)}$.

\ep

The following result is obvious.
  \ble\lb{d}
 Let $U$ be a $\phi_0$-invariant subspace of $V_0$ and $\mathcal{B}=\{v_1,\cdots,v_k\}$ be a basis of $U$. Assume $M(\phi_0|_U, \mathcal{B})=A$. Then  $J(U)$ is also a $\phi_0$-invariant subspace of $V_0$, with basis $J(\mathcal{B})=\{Jv_1,\cdots,Jv_k\}$, and $M(\phi_0|_{J(U)}, J(\mathcal{B}))=\bar{A}$.
 In particular, if $M(\phi_0|_U, \mathcal{B})=\mathrm{J_m(\la)}$ for some $\la\in \mathbb{C}$, then $M(\phi_0|_{J(U)}, J(\mathcal{B}))=\mathrm{J_m(\bar{\la})}$.
 \ele
Now we can prove the Jordan canonical form theorem for quaternionic linear operators.
\bthm\lb{r}
(1) Let $V$ be a quaternionic linear space and $\phi\in End_\mathbb{H}(V)$. Then there exist $\phi$-invariant subspaces $V_i$ of $V$, $i=1,\cdots,k$, such that $V=\oplus_{i=1}^k V_i$ and there exists a basis $\mathcal{B}_i$ for $V_i$ such that $M(\phi|_{V_i},\mathcal{B}_i)$ is a Jordan block $\mathrm{J_{m_i}(\la_i)}$, where $\la_i\in \mathbb{H}$.

(2)Such Jordan blocks are uniquely determined by $\phi$ up to permutation and up to a replacement of $\la_1,\cdots,\la_k$ with $p_1^{-1}\la_1 p_1,\cdots,p_k^{-1}\la_k p_k$ within the blocks $\mathrm{J_{m_1}(\la_1)},\cdots , \mathrm{J_{m_k}(\la_k)}$ respectively, where $\la_i\in \mathbb{H}$ and $p_i\in \mathbb{H}\setminus \{0\}$.
\ethm

\bp
(1)By last theorem, $V_0$ can be written as a direct sum of $J$-paired Jordan subspaces, i.e., $V_0=\bigoplus_{i=1}^k (W_i\oplus JW_i)$, where $W_i$ are Jordan subspaces of $\phi$. Let $\mathcal{B}_i$ be a Jordan chain in $W_i$. Then $J(\mathcal{B}_i)$ is a Jordan chain in $J(W_i)$. Let the basis $\mathcal{B}_0$ of $V_0$ be the union of all the  Jordan chains $\mathcal{B}_i$ and $J(\mathcal{B}_i)$. Then
 $M(\phi_0,\mathcal{B}_0)=\bigoplus_{i=1}^k (\mathrm{J_{m_i}(\la_i)}\oplus \mathrm{J_{m_i}(\bar{\la_i}))}$, where $\bigoplus_{i=1}^k (\mathrm{J_{m_i}(\la_i)}\oplus \mathrm{J_{m_i}(\bar{\la_i}))}$ denotes the block diagonal matrix $diag(\mathrm{J_{m_1}(\la_1)},\mathrm{ J_{m_1}(\bar{\la_1})},\cdots,\mathrm{J_{m_k}(\la_k)}, \mathrm{J_{m_k}(\bar{\la_k}))}$, is a Jordan decomposition for $\phi_0$. Let $V_i=W_i\oplus J(W_i),$ then $V_i$ is a $\mathbb{H}$-subspace of $V$ and $V=\bigoplus_{i=1}^k V_i$. Then  the union $\mathcal{B}$ of all the  Jordan chains $\mathcal{B}_i$, for $i=1,\cdots,k$, is a basis for $V$ and $M(\phi,\mathcal{B})=\bigoplus_{i=1}^k \mathrm{J_{m_i}(\la_i)}$ is a desired Jordan decomposition for $\phi$.

(2)  Assume  $V=\oplus_{i=1}^s V^{'}_i$ is another such decomposition of $V$ into $\phi$-invariant Jordan subspaces $V^{'}_i$, and $\mathcal{B}^{'}_i$ is the corresponding Jordan chain, for $i=1,\cdots, s$. Assume $M(\phi|_{V^{'}_i},\mathcal{B}^{'}_i)=\mathrm{J_{l_i}(\mu_i)}$. Let $\mathcal{B}^{'}$ be the union of all the  Jordan chains $\mathcal{B}^{'}_i$. Then $M(\phi,\mathcal{B}^{'})=\bigoplus_{i=1}^s \mathrm{J_{l_i}(\mu_i)}$.

Assume  $\mathcal{B}^{'}_i=\{v_1,\cdots, v_{l_i}\}$ and $M(\phi,\mathcal{B}^{'}_i)=\mathrm{J_{l_i}(\mu_i)}$. One has $p_i^{-1}\mu_i p_i=\nu_i\in \mathbb{C}$ for suitable choices of nonzero $p_i\in \mathbb{H}$. Let $\mathcal{B}^{''}_i=( v_1 p_i,\cdots, v_{l_i} p_i)$.  then $M(\phi|_{V^{'}_i},\mathcal{B}^{''}_i)=\mathrm{J_{l_i}(\nu_i)}$ by (\ref{s}). Let $\mathcal{B}^{''}$ be the union of $\mathcal{B}^{''}_i$, which is a basis for $V$, then $M(\phi,\mathcal{B}^{''})=\bigoplus_{i=1}^s \mathrm{J_{l_i}(\nu_i)}$ is also a Jordan decomposition for $\phi$. Let the basis $\mathcal{B}^{''}_0$ of $V_0$ be the union of all the  Jordan chains $\mathcal{B}^{''}_i$ and $J(\mathcal{B}^{''}_i)$. Then $M(\phi_0, \mathcal{B}^{''}_0 )=\bigoplus_{i=1}^s (\mathrm{J_{l_i}(\nu_i)}\oplus \mathrm{J_{l_i}(\bar{\nu_i}))}$ is a Jordan decomposition for $\phi_0$.

By (1),  $M(\phi_0,\mathcal{B}_0)=\bigoplus_{i=1}^k (\mathrm{J_{m_i}(\la_i)}\oplus \mathrm{J_{m_i}(\bar{\la_i}))}$ is also a Jordan decomposition for $\phi_0$. By the uniqueness of the Jordan blocks of the complex linear operator $\phi_0$, one has $s=k$ and there is a permutation $\sigma$ on $\{1,\cdots,k\}$ such that ${l_i}=m_{\sigma i}$ and $\nu_i=\la_{\sigma i}$ or $\nu_i=\overline{\la_{\sigma i}}$ for $i=1,\cdots,k$. As $\la_{\sigma i}$ and $\overline{\la_{\sigma i}}$ are similar in $\mathbb{H}$, one can choose $p_i$ properly such that $\nu_i=\la_{\sigma i}$. This completes the proof.
\ep

\bco The exponential maps $exp_1:M_n(\mathbb{H})\rt GL_n(\mathbb{H})$ and $exp_2:End_\mathbb{H}(V)\rt GL(V)$ are both surjective.\eco
\bp
We will show that $exp_1:M_n(\mathbb{H})\rt GL_n(\mathbb{H})$  is surjective. Then it follows that $exp_2:End_\mathbb{H}(V)\rt GL(V)$ is also surjective by Figure 1.

As $exp_1(PAP^{-1})=P exp_1(A)P^{-1}$ for any $P\in GL_n(\mathbb{H})$ and $A\in M_n(\mathbb{H})$, one only needs to show that for any $A\in GL_n(\mathbb{H})$ with $A= \oplus_{i=1}^k J_{m_i}(\la_i)$ has a preimage under $exp_1$, by the above theorem. Moreover, one can assume that $\la_i\in \mathbb{C}$ for each $i$ by Theorem \ref{r} (2). Then it suffices to show that each Jordan block $J_{m_i}(\la_i)$ has a preimage under $exp_1$. But this follows from the surjectivity of the exponential map $exp:M_{m_i}(\mathbb{C})\rt GL_{m_i}(\mathbb{C})$.
\ep

\end{document}